\documentclass{proc-l}
\usepackage{amssymb}
\usepackage{eufrak}
\usepackage{amsmath}
\usepackage{shadow}
\begin{document}
\newtheorem{theorem}{Theorem}[section]
\newtheorem{lemma}[theorem]{Lemma}
\newtheorem{cor}[theorem]{Corollary}
\theoremstyle{definition}
\newtheorem{definition}[theorem]{Definition}
\newtheorem{example}[theorem]{Example}
\newtheorem{xca}[theorem]{Exercise}

\theoremstyle{remark}
\newtheorem{remark}[theorem]{Remark}

\numberwithin{equation}{section}

\def\CX{{\mathbb C}}
\def\QX{{\mathbb Q}}
\def\NX{{\mathbb N}}
\def\ZX{{\mathbb Z}}
\def\DX{{\mathbb D}}
\def\AX{{\mathbb A}}
\def\ord{\mbox{ord }}
\def\GL{{\rm GL}}
\def\SL{{\rm SL}}
\def\gl{{\rm gl}}
\def\sl{{\rm sl}}
\def\M{{\rm gl}}
\def\d{{
\partial}}
\def\ov{\overline}
\def\frakk{{\mathfrak k}}
\def\frakK{{\mathfrak K}}
\def\frakE{{\mathfrak E}}
\def\frakF{{\mathfrak F}}
\def\frakU{{\mathfrak U}}
\def\ofrakF{{\overline{\mathfrak F}}}
\def\oofrakF{{\overline{\mathfrak F}}_0}
\def\ofrakG{{\overline{\mathfrak G}}}
\def\oofrakG{{\overline{\mathfrak G}}_0}
\def\frakG{{\mathfrak G}}
\def\td{{\tilde{D}}}
\def\ttd{{\tilde{\tilde{D}}}}
\def \ldf{{{\rm LDF}_0}}
\def\QED{\hbox{\hskip 1pt \vrule width4pt height 6pt depth 1.5pt \hskip 1pt}}
\def\calD{{\mathcal D}}
\def\calA{{\mathcal A}}
\def\calC{{\mathcal C}}
\def\calS{{\mathcal S}}
\def\calT{{\mathcal T}}
\def\calF{{\mathcal F}}
\def\calE{{\mathcal E}}
\def\calL{{\mathcal L}}
\def\calG{{\mathcal G}}
\def\calK{{\mathcal K}}
\def\calO{{\mathcal O}}
\def\calU{{\mathcal U}}
\def\calM{{\mathcal M}}
\def\calP{{\mathcal P}}
\def\calR{{\mathcal R}}
\def\calL{{\mathcal L}}
\def\calB{{\mathcal B}}
\def\calN{{\mathcal N}}
\def\calCF{{\calC_0(\calF)}}
\def\gal{{\rm DiffGal}}
\def\C{{\rm Const}}
\def\d{{
\partial}}
\def\dx{{
\partial_{x}}}
\def\dt{{
\partial_{t}}}
\def\Autd{{{\rm Aut}_\Delta}}
\def\Gal{{\rm Gal}}
\def\PGal{{\rm Gal_\Delta}}
\def\Ga{{\bf G_a}}
\def\Gm{{\bf G_m}}
\def\ld{{l\d}}
\def\vr{{\Vec{r}}}
\def\vta{{\Vec{\tau}}}
\def\vt{{\Vec{t}}}
\def\KPV{{K^{\rm PV}_A}}
\def\ord{{\rm ord}}
\def\kbar{{\bar{k}}}

\title[Model Theory of Partial Differential Fields]{Model Theory of Partial Differential Fields:\\ From Commuting to Noncommuting Derivations}
\author{Michael F.
Singer}
\address{North Carolina State University\\ Department of
Mathematics\\Box 8205\\ Raleigh, North Carolina 27695-8205}
\email{ 
singer@math.ncsu.edu} \thanks{The preparation of this article was partially supported by NSF Grant CCR-0096842.}
\subjclass[2000]{Primary: 03C10; Secondary: 35A05, 12H05}
\date{}
\commby{Julia F. Knight}
\begin{abstract} McGrail \cite{mcgrail} has shown the existence of a model completion for the universal theory of fields on which a finite number of commuting derivations act and, independently, Yaffe \cite{yaffe} has shown the existence of a model completion  for the univeral theory of  fields on which a fixed Lie algebra acts as derivations. We show how to derive the second result from the first. 
\end{abstract}
\maketitle
\section{Introduction.} In \cite{mcgrail}, McGrail gives axioms for the model completion of the universal theory of fields of characteristic zero with several commuting derivations. Independently, Yaffe \cite{yaffe} gave  axioms for the model completion of the universal theory of a more general class of fields: $LDF_0$, the universal theory of fields of characteristic zero together with a finite dimensional Lie algebra acting as derivations.  The goal of this short note is to show that starting from McGrail's results one can quickly write down axioms for Yaffe's model completion, that is, one can reduce the noncommutative case to the commutative case.

Another axiomatization of the model completion of a theory of differential fields with noncommuting derivations   has been given  by Pierce \cite{pierce} using differential forms and a version of the Frobenius Theorem. 
In fact, one can state the Frobenius Theorem as a result that allows one to replace noncommuting vector fields with commuting ones: {\em given an involutive analytic system $\calL_y, y\in M$ of tangent spaces of rank p on an analytic manifold M and a point $x\in M$, one can choose analytic coordinates $x_1, \ldots ,x_p$ around $x$ such that $\calL_y$ is spanned by $(\d/\d x_1)_y, \ldots , (\d/\d x_p)_y$ for all $y$ in an open set containing $x$} ({\em cf.} \cite{varadarajan_lie}, Theorem 1.3.3).
In his introductory remarks, Pierce states  that one might be able to reduce the noncommutative case to the commutative case (\cite{pierce}, pp.924-925) but does not give details, proceeding rather to develop his results {\em ab initio}.  The present note carries through this reduction.

I would like to thank Evelyne Hubert for several useful discussions and especially for pointing out  the relation between the  result of Cassidy and Kolchin appearing as Proposition 6, p. 12 of \cite{kolchin_groups} and  the Frobenius Theorem. I would also like to thank Yoav Yaffe for discussions that helped clarify my initially vague ideas and the anonymous referee for  useful comments.

\section{Commuting Bases of Derivations} In \cite{yaffe}, Yaffe defines the universal theory $LDF_0$ of Lie differential fields. To do this he fixes a field $\calF$ of characteristic zero, a finite dimensional $\calF$-vector space $\calL$ with a Lie multiplication making it a Lie algebra over a subfield of $\calF$ and a vector space homomorphism $\phi_\calF:\calL \rightarrow {\rm Der}(\calF)$, the Lie algebra of derivations on $\calF$,  preserving the  Lie multiplication.  We fix a basis $\{D_1, \ldots , D_n\}$ of $\calL$ and let
$[D_k,D_l] = \sum_m \alpha_{kl}^mD_m$ for some $\alpha_{kl}^m\in \calF$.  The elements $\alpha_{kl}^m$ are called the structure constants of $\calL$. The language for $LDF_0$ is the language of rings together with unary function symbols $D_i, i = 1, \ldots , n$ and constant symbols for each element of $\calF$. The theory $LDF_0$ consists of
\begin{itemize}
\item the diagram of $\calF$ including the action of the $D_i$,
\item the theory of integral domains of characteristic zero,
\item axioms stating that the $D_i$ are derivations,
\item for each $k,l$ an axiom of the form $\forall x(D_kD_lx-D_lD_kx = \sum_m \alpha_{kl}^mD_m x)$.
\end{itemize}
We shall refer to a model of $LDF_0$ as  a {\em Lie ring over $\calL$} or, more simply, a {\em Lie ring}, if it is clear which $\calL$ is used. 

We will need another concept from \cite{yaffe}: the ring of normal polynomials. This is defined as follows. Let $\calA$ be a model of $LDF_0$. For each $I = \langle i_1, \ldots ,i_n\rangle \in \NX^n$ define a variable $X_I$ and consider the ring $\calR_\calA = \calA[\{X_I\}_{I\in \NX^n}]$ of polynomials in this infinite set of variables. Yaffe shows (\cite{yaffe}, p. 57-61) how one can define an action of $\calL$ on $\calR_\calA$ so that this ring is a model 
of $LDF_0$ and one sees from his proof that the $D_i$ act as linearly independent derivations. Although this is not needed in the following lemma, we note that Yaffe also shows  that if $\calB$ is a model of $LDF_0$ that extends $\calA$ and $b \in \calB$, then the map from $\calR_\calA$ to $\calB$ given by $X_I \mapsto D^I(b) = D_1^{i_1}D_2^{i_2}\ldots D_n^{i_n}(b)$ for all $I$, is a homomorphism of structures. Yaffe's construction of $\calR_\calA$ allows us to  conclude
\begin{lemma}\label{lem1}If $\calA$ is a  model of
$LDF_0$ with derivations $D^\calA_1, \ldots ,D^\calA_n$,
then
$\calA$ embeds in a model $\calB$ of $LDF_0$ with derivations $D_1^\calB,
\ldots D_n^\calB$ extending the derivations of $\calA$ such that the 
$D_1^\calB, \ldots D_n^\calB$ are linearly independent over
$\calB$.
\end{lemma}

To show that the axioms in the next section yield a model completion of this 
theory we need two additional facts that are contained in the next results.

\begin{lemma} \label{lem2} Let $\calA$ with derivations $D^\calA_1,
\ldots , D^\calA_n$ be a model of $LDF_0$ and assume that the $D_i^\calA$ are
linearly independent over $\calA$. Then there exist $a_{i,j} \in \calA$ such that
the derivations $\overline{D}_i^\calA = \sum_j
a_{i,j} D_j^\calA$ are linearly independent over $\calA$ and commute.
\end{lemma}
\begin{proof} We shall assume $\calA$ is a field and follow the proof of Proposition
6 in Chapter 0, \S5 of \cite{kolchin_groups} with a few small modifications.  As noted in the introduction, this is an algebraic version of the Frobenius Theorem. 

Let $\calC = \cap_{i=1}^n {\rm Ker }\ D^\calA_i$.  Let 
$\{x_i\}_{i \in I}$ be a transcendence basis of 
$\calA$ over $\calC$.  Let $\delta_i$ be the unique derivation on
$\calA$ satisfying $\delta_i(c) = 0$ for all $c\in \calC$ and
$\delta_i(x_j) = 1$ if $i=j$ and $0$ if $i\neq j$. Any
$x \in \calA$ will lie in a subfield algebraic over
$\calC(x_1, \ldots x_N) $ for some $N$ and so $\delta_i(x) = 0$
for all $i>N$. For any
derivation $D$ of $\calA$ that is trivial on $\calC$, consider the
sum
$\sum_{i \in I}D(x_i) \delta_i$.  Although this sum is
infinite, the previous remark shows that for any $x \in \calA,$
the sum $\sum_{i \in I}D(x_i) \delta_i$ makes sense and using the results of Chapter VII, $\S5$
of \cite{LANG}, one can show that 
$D = \sum_{i \in I}D(x_i) \delta_i$.

Since the $D_i$ are linearly independent over $\calA$ there exist 
$a_{i,j} \in \calA$ such that
the derivations $\ov{D}_i^\calA = \sum_j
a_{i,j} D_j^\calA$ satisfy (after a possible renumbering of the $\delta_j$)
\begin{equation*}\ov{D}_i^\calA = \delta_i \mbox{ (mod } \sum_{j>n} \calA\delta_j) \end{equation*}
for $i= 1, \ldots ,n$. To see that the  $\ov{D}^\calA_i$ are linearly independent
over $\calA$ note that if  $\sum_j b_j \ov{D}_j^\calA$ = 0 then $0 = \sum_j b_j \ov{D}_j^\calA(x_i) =
b_i$.  I now claim that the $\ov{D}^\calA_i$ commute.  Since the $\delta_i$ commute, we have that
\begin{equation*}[\ov{D}_r^\calA,\ov{D}_s^\calA] = 0 \mbox{ (mod }\sum_{i>n}\calA\delta_i) \end{equation*}
for any $r,s$. Since the
$\calA$-span of the $\ov{D}^\calA_i$ is closed under $[ \ , \ ]$, we have that 
there exist $b_{r,s,j} \in \calA$ such that
\begin{eqnarray*} [\ov{D}_r^\calA,\ov{D}_s^\calA] & = & \sum_{j=1}^n b_{r,s,j}\ov{D}^\calA_j\\
   & = & \sum_{j=1}^n b_{r,s,j}\delta_j\mbox{ (mod
}\sum_{i>n}\calA\delta_i) \ \ .
  \end{eqnarray*}
Therefore, we have that $b_{r,s,j} = 0$ for all $r,s,$ and
$j$ and hence that the $\ov{D}^\calA_j$ commute. 
\end{proof}
We shall need to compare Lie rings for two different Lie algebras $\calL_1$ and $\calL_2$.  I will denote by
$LDF_0^1$ (resp. $LDF_0^2$) the theory of Lie rings based on the
action of the Lie algebra $\calL_1$ (resp. $\calL_2$). I will
assume the two  algebras are of the same dimension over $\calF$. 
\begin{lemma}\label{lem3} Let $\calA$ with derivations $D^1_1,
\ldots , D^1_n$ be a model of $LDF_0^1$ and assume that the $D_i^1$ are
linearly independent over $\calA$.  For $i = 1, \ldots, n$, let $D_i^2 = \sum_j
a_{i,j} D_j^1$ for some $a_{i,j} \in \calA$ and assume that $\calA$ with
derivations  $D^2_1,
\ldots , D^2_n$ is a model of $LDF_0^2$. If $\calB$ is an extension of $\calA$ that
is a model of $LDF_0^1$ with respect to the extensions of the derivations
$D^1_1, \ldots , D^1_n$, then the formulas $D_i^2 = \sum_j
a_{i,j} D_j^1$ define derivations on $\calB$ such that $\calB$ with these
derivations is a model of $LDF_0^2$ 
\end{lemma}
\begin{proof}Let $\alpha_{k,l}^m$ be the structure constants of
$\calL_1$ and $\beta_{k,l}^m$ be the structure constants of $\calL_2$. We
have
\begin{eqnarray*}
[D^2_l,D^2_k] &=& \sum_j(\sum_i(a_{l,i}D^1_i(a_{k,j})-a_{k,i}D^1_i(a_{l,j})) + \sum_{r,s}
a_{l,r}a_{k,s}\alpha_{r,s}^j)D^1_j, \mbox{ and } \\
\sum_m \beta_{l,k}^m D^2_m & = & \sum_j (\sum_m \beta_{l,k}^m a_{m,j})D^1_j
\end{eqnarray*}
First note that since the  $D^1_j$ are linearly independent over $\calA$, these derivations are linearly independent over $\calB$ as well. We therefore have that $[D^2_l,D^2_k] = \sum_m \beta_{l,k}^m D^2_m$ if and
only if 
\begin{equation*}\sum_i(a_{l,i}D^1_i(a_{k,j})-a_{k,i}D^1_i(a_{l,j}) )+ \sum_{r,s}
a_{l,r}a_{k,s}\alpha_{r,s}^j= \sum_m \beta_{l,k}^m a_{m,j} \end{equation*}
for all $j$.  If this holds in $\calA$ then it will hold in $\calB$. 
\end{proof}
In particular, if the basis elements  $D_i^2$  commute, we have
\begin{cor} \label{cor} Let $\calA$ with derivations $D^1_1,
\ldots , D^1_n$ be a model of $LDF_0$ and assume that the $D_i^1$ are
linearly independent over $\calA$.  For $i = 1, \ldots, n$, let $D_i^2 = \sum_j
a_{i,j} D_j^1$ for some $a_{i,j} \in \calA$ and assume that the $D_i^2$ commute as derivations on $\calA$. If $\calB$ is an extension of $\calA$ that
is a model of $LDF_0$, then the $D_i^2$ commute as derivations on $\calB$.
\end{cor} 
Note that in model theoretic terms,  the statement that the basis elements $D_i^2$ commute, which {\em a priori} is a universal statement, is actually (equivalent to) a quantifier free statement, given that the  $D_i^1$ are linearly independent.
\section{Axioms for the Model Completion of $LDF_0$} Let us begin by recalling the situation for fields with commuting derivations. In \cite{mcgrail}, McGrail gave axioms for the model completion
$m$-DCF of the universal theory $m$-DF of differential fields with
$m$ commuting derivations $D_1, \ldots ,D_n$. She also showed that this former
theory has elimination of quantifiers.  In particular, for any
field
$F$ with $m$ commuting derivations and any system $\calS(u_1,
\ldots, u_r, v_1, \ldots ,v_s) = \{f_1(u_1,
\ldots, u_r, v_1, \ldots ,v_s) = 0, 
\ldots , f_t(u_1, \ldots, u_r, v_1, \ldots ,v_s) = 0, g(u_1,
\ldots, u_r, v_1, \linebreak \ldots ,v_s) \neq 0\}$ of differential
polynomials  there
exist systems 
$\calT_1(u_1, \ldots, u_r), \ldots , \linebreak \calT_l(u_1,
\ldots, u_r)$ such that for any  $u_1, \ldots u_r \in F$ there
exists a solution $v_1, \ldots , v_s$ in some differential
extension of $F$  if and only if the
$u_i$ satisfy one of the systems $\calT_i(u_1, \ldots,
u_r)$.

Let $\calF, \calL$ and $\phi_{\calF}$ be as above.  Let $\calR_\calF$ be the ring of normal polynomials with coefficients in $\calF$. For any $t \in \NX$, we will denote by   $\calR_\calF[x_1, \ldots , x_t]$ the usual (not differential) ring of polynomials in the variables $x_1, \ldots, x_t$ with coefficients in $\calR_\calF$. The axioms for
$LDCF_0$, the theory of Lie differentially closed fields of characteristic zero,  are  the axioms for
$LDF_0$ plus the axioms for fields and the  following axioms:
\begin{enumerate}
\item for any $t$ and any polynomial $p(x_1, \ldots ,x_t, \{X_I\}) \in \calR_\calF[x_1, \ldots , x_t]$ we have an axiom that
states that for any $a_1, \ldots , a_n$ such that $p(a_1, \ldots ,a_n, \{X_I\})
\neq 0$ there exists  $b$ such that $p(a_1,$ $ \ldots ,a_n, \{D^I(b)\}) \neq 0$ 
\item there exists an $n^2$-tuple of elements  ${\bf x} =
(x_{1,1},
\ldots , x_{n,n})$ such that the derivations $D_1^{\bf{x}} =
\sum_j x_{1,j}D_j, \ldots , D_n^{\bf{x}} =
\sum_j x_{n,j}D_j$ form a linearly independent set of commuting 
derivations.

\item for any ${\bf x}$  such that $\{D_1^{\bf{x}},  \ldots
, D_n^{\bf{x}}\}$ is a linearly independent  set of commuting derivations and
any system  $S(u_1, \ldots ,u_n, v_1,
\ldots ,v_m)$ involving differential polynomials
in the $D_i^{\bf x}$,  we have  
\begin{equation*}\exists u_1,
\ldots v_m
\calS(u_1,
\ldots ,u_n, v_1,
\ldots ,v_m)\end{equation*}  if and only if \begin{equation*}\bigvee_{j=1}^l \calT_j(u_1,
\ldots ,u_n)\end{equation*}
where the $\calT_j$ are as described above.
\end{enumerate}
Note that (3) implies that the models of this theory are algebraically
closed.  

I will use Blum's criteria (\cite{sacks}, Section 17) to show that these axioms give the model
completion of $LDF_0$.  The first step is to show that any model $\calA$ of
$LDF$ can be extended to a model of $LDCF_0$. Taking quotient fields of the ring constructed in Lemma~\ref{lem1} we
extend $\calA$ to a model $K^0$ of $LDF_0$ were the derivations are linearly
independent over $K^0$. Note that the ring constructed in Lemma~\ref{lem1} allows us to conclude that $K^0$  satisfies
axiom scheme (1).  Lemma~\ref{lem2} implies that axiom (2)  holds. Let
$\calN_0$ be the set of $n^2$-tuples of elements in $K^0$ and assume that this set is
well ordered. Let ${\bf x}$ be the smallest element of $\calN_0$
 such that
$\{D_1^{\bf{x}}, 
\ldots , D_n^{\bf{x}}\}$ is a linearly independent  set of commuting derivations
and let
$K^0_1$ be the differential closure of $K^0$ thought of as a field with commuting
derivations $\{D_1^{\bf{x}},  \ldots
, D_n^{\bf{x}}\}$. Since the $D_i$ can be expressed a $K^0$-linear combinations
of the $D_i^{\bf x}$, Lemma~\ref{lem3} implies that $K^0_1$ is still  a model of $LDF_0$. 
Let ${\bf \bar x}$ be the next smallest element of $\calN_0$
 such that
$\{D_1^{\bf{\bar x}}, 
\ldots , D_n^{\bf{\bar x}}\}$ is a linearly independent  set of commuting derivations
and let
$K^0_2$ be the differential closure of $K^0_1$ thought of as a field with commuting
derivations $\{D_1^{\bf{ \bar x}},  \ldots
, D_n^{\bf{ \bar x}}\}$. We can continue in this way for all elements of $\calN_0$ and,
taking unions, form a field $\bar{K}^0$ such that for any ${\bf x}$  in $\calN_0$
 such that
$\{D_1^{\bf{x}}, 
\ldots , D_n^{\bf{x}}\}$ is a linearly independent  set of commuting derivations,
$\bar{K}^0$ contains the differential closure of $K^0$  thought of as a field with commuting
derivations $\{D_1^{\bf{x}},  \ldots
, D_n^{\bf{x}}\}$. Again Lemma~\ref{lem3} implies that $\bar{K}^0$ is a model of $LDF_0$. We let $K^1$ be the quotient field of $\calR_{\bar{K}^0}$.  Note that $K^1$ satisfies axiom scheme (1).  We now repeat
this process  and form a field
$K^2$ such  that for any $n^2$-tuple ${\bf x}$  of elements in $K^1$
 such that
$\{D_1^{\bf{x}}, 
\ldots , D_n^{\bf{x}}\}$ is a linearly independent  set of commuting derivations,
$K^2$ contains the differential closure of $K^1$  thought of as a field with commuting
derivations $\{D_1^{\bf{x}},  \ldots
, D_n^{\bf{x}}\}$ and also satisfies axiom scheme (1). One now sees  that $K^{\infty} = \cup K^i$ is a model of $LDF_0$
and satisfies the axioms of $LCDF_0$.

I will now show that if $\calA \models LDF_0$, $\calB \models LDCF_0$, $\calB$ being
$|\calA|^+$-saturated and
$\calA
\subset
\calB$, then for any simple extension $\calA(a)$ there is an $\calA$-embedding
$f:\calA(a) \rightarrow \calB$.  I will first show that I can assume that $\calA$
satisfies axioms (1) and (2) above. 

Let $\calR_{\calA(a)} = \calA(a)(\{X_I\})$ be the  field fractions of normal polynomials over $\calA(a)$.  Axioms (1) and saturation imply that there
is an $\calA$-embedding of $\calA(\{X_I\})$ into $\calB$. Furthermore, $\calA(\{X_I\})$ satisfies axiom schemes (1) and  (2).  If we can extend
this embedding to 
$\calA(\{X_I\})(a)$, then restricting to
$\calA(a)$ will give the desired conclusion.  We therefore will assume from the
beginning that
$\calA$ satisfies axiom schemes (1) and (2).

Let ${\bf x}$ be elements of $\calA$ such that $D_i^{\bf x}$ are linearly independent
commuting derivations. Corollary~\ref{cor} implies that these derivations commute on $\calA(a)$ as well. The isomorphism type of $a$ over the field $\calA$ with
derivations
$\{D^{\bf x}_i\}$ is determined by the set of $\{D^{\bf x}_i\}$-differential
polynomials that $a$ satisfies and the set that it does not satisfy. By axiom scheme (3) and
saturation, we can realize this type in $\calB$ and so get an embedding of
$\calA(a)$ into $\calB$, considered as  $\{D^{\bf x}_i\}$-differential fields.  Since
the $\{D_i\}$ are linear combinations of the $\{D^{\bf x}_i\}$, this is also an
embedding as models of $LDF_0$. 
\newcommand{\SortNoop}[1]{}
\providecommand{\bysame}{\leavevmode\hbox to3em{\hrulefill}\thinspace}
\providecommand{\MR}{\relax\ifhmode\unskip\space\fi MR }
\providecommand{\MRhref}[2]{%
  \href{http://www.ams.org/mathscinet-getitem?mr=#1}{#2}
}
\providecommand{\href}[2]{#2}

\end{document}